# Galerkin-Bernstein Approximations for the System of Third-Order Nonlinear Boundary Value Problems

**Snigdha Dhar and Md. Shafiqul Islam***

Department of Applied Mathematics, University of Dhaka, Dhaka-1000, Bangladesh

*Corresponding Author:* mdshafiqul@du.ac.bd

**Abstract**

This paper is devoted to find the numerical solutions of one dimensional general nonlinear system of third-order boundary value problems (BVPs) for the pair of functions using Galerkin weighted residual method. We derive mathematical formulations in matrix form, in details, by exploiting Bernstein polynomials as basis functions. A reasonable accuracy is found when the proposed method is used on few examples. At the end of the study, a comparison is made between the approximate and exact solutions, and also with the solutions of the existing methods. Our results converge monotonically to the exact solutions. In addition, we show that the the derived formulations may be applicable by reducing higher order complicated BVP into a lower order system of BVPs, and the performance of the numerical solutions is satisfactory.

## 1 Introduction

Ordinary differential systems have received a lot of interest in studies as a result of their frequent occurrence in numerous applications. Second-order ordinary differential systems can simulate a number of natural phenomena. For instance, while studying chemically reacting systems, fluid mechanics, relativistic mechanics, gas dynamics, and nuclear physics, the known Emden-Fowler equations appear. Numerous authors have made an effort to increase accuracy in the literature on numerical analysis. Linear second-order boundary value problems can be effectively solved using the finite difference approach. However, when it comes to deal with nonlinear equations, solving the corresponding boundary value problems can be quite challenging. The Galerkin method [1] have been used to solve the two point BVP by the authors Bhatti and Bracken [2]. Bernoulli polynomials [3] have been used to solve second-order both linear and nonlinear BVPs with Dirichlet, Neumann and Robin boundary conditions. Spline functions [4, 5] have been extensively researched because piecewise polynomials [6] can be differentiated, integrated, and approximate to any function with any desired accuracy.

Second order ordinary differential systems are used to describe a variety of problems in biology, engineering, and physics. A lot of works have been done for obtaining numerical solutions for the linear system of second-order boundary value problems such as, Geng and Cui [7] represented the approximate solutions of system of second-order linear and nonlinear BVP in the form of series in the reproducing kernel space. For solving a nonlinear system of second-order BVP, the variational iteration technique [8] was introduced. A nonlinear system of second order BVPs was solved by Dehghan and Saadatmandi using the sinc-collocation method [9]. The exact and approximate solutions were expressed in the reproducing kernel space by Du and Cui [10]. A novel strategy to solve the nonlinear systems of second order BVPs is provided, depending on the homotopy perturbation method and the reproducing kernel method. The advantages of both of these techniques are combined in the homotopy perturbation-reproducing kernel method [11], which can be utilised to effectively handle systems of nonlinear boundary value problems. A few examples of linear and nonlinear systems were solved using the Galerkin weighted residual technique [6]. In order to solve systems of second-order BVPs, a family of boundary value methods was used in [12], as a block unification approach. In order to solve systems of singular boundary value problems, an optimization approach [13] was presented and then solved via continuous genetic algorithm [14].

Recently, Galerkin finite element method [15] has been used to determine the approximate solutions of a coupled second-order BVPs. In recent years, coupled Lane-Emden [16, 17] boundary value problems have received a lot of attention. The singularity is the primary barrier to solutions, and many authors are working to solve it. In order to efficiently solve the system of Lane-Emden type equations that arise in physics, star structure, and astrophysics, Ala'yed, Saadeh and Qazza [19] proposed a collocation approach based on cubic B-spline functions.

Solving a nonlinear system of third-order boundary value problems is quite difficult. Consequently, a few researchers have attempted to solve numerically. Ezzati and Aqhamohamadi [21] used He's Homotopy Perturbation method to solve the nonlinear system of third-order boundary value problems. Block method [22] was used to solve the nonlinear system of third-order boundary value problems. We observe that a little attention has been given for solving third-order system of boundary value problems. Thus we are motivated to solve the system of third-order boundary value problems by Modified Galerkin weighted residual method with Bernstein polynomials as trial functions. However, we organize this research work as following.

First of all, we give basic idea on Bernstein polynomials in section 2. Mathematical formulations are described elaborately, in matrix form, in section 3. Section 4 is devoted to numerical experiments and the discussion of the obtained results in tabular form and graphically. Then in section 5, we provide an application that sixth order BVP may be modeled into a system of third order BVPs which can be solved by the proposed technique. Finally, conclusions and references are amended.

## 2  Bernstein Polynomials

The general form of the Bernstein polynomials of degree $n$ over the interval $[a,b]$ is defined by [3, 6, 23]

$$\phi_{i,n}(x) = \binom{n}{i} \frac{(x-a)^i (b-x)^{n-i}}{(b-a)^n}, \ a \leq x \leq b \quad i = 0, 1, 2, \ldots, n$$

It is important to observe that each of these $n+1$ polynomials, with a degree of $n$, fulfills the following properties:

(i)  $\phi_{i,n}(x) = 0$ \qquad if $i < 0$ or $i > n$

(ii)  $\sum_{i=0}^{n} \phi_{i,n}(x) = 1$

(iii)  $\phi_{i,n}(a) = \phi_{i,n}(b) = 0$ \qquad $1 \leq i < n$

For simplification, we denote $\phi_{i,n}(x)$ by $\phi_i(x)$. The Bernstein polynomials of degree 3, 4 and 5 are given below, respectively

**Degree 3:**    $\phi_0(x) = (1-x)^3$    $\phi_1(x) = 3x(1-x)^2$    $\phi_2(x) = 3x^2(1-x)$    $\phi_3(x) = x^3$

**Degree 4:**    $\phi_0(x) = (1-x)^4$    $\phi_1(x) = 4x(1-x)^3$    $\phi_2(x) = 6x^2(1-x)^2$    $\phi_3(x) = 4x^3(1-x)$

             $\phi_4(x) = x^4$

**Degree 5:**    $\phi_0(x) = (1-x)^5$    $\phi_1(x) = 5x(1-x)^4$    $\phi_2(x) = 10x^2(1-x)^3$    $\phi_3(x) = 10x^3(1-x)^2$

             $\phi_4(x) = 5x^4(1-x)$    $\phi_5(x) = x^5$

To solve a BVP using the Galerkin method, it is required that each of these polynomials satisfies the homogeneous representation of the essential boundary conditions. Here, $\phi_0(x)$ and $\phi_n(x)$ do not satisfy the homogeneous boundary conditions. Therefore, in order to satisfy the homogeneous boundary conditions, we use only $\phi_i(x)$ for $1 \leq i \leq (n-1)$. We use Bernstein polynomials of degree 3, 4 and 5 throughout this paper.

## 3   Mathematical Formulation of System of Third-Order BVPs

In recent times, the interest in boundary value problems of a system of ordinary differential equations has been sparked among researchers in mathematics, physics, engineering, biology, and other fields.

The general linear system of two third order differential equations in two unknown functions $p(x)$ and $q(x)$ of the form below is taking into consideration [6]

$$\begin{cases} a_1(x)p''' + a_2(x)p'' + a_3(x)p' + a_4(x)p + a_5(x)q''' + a_6(x)q'' + a_7(x)q' + a_8(x)q = g_1(x) \\ b_1(x)p''' + b_2(x)p'' + b_3(x)p' + b_4(x)p + b_5(x)q''' + b_6(x)q'' + b_7(x)q' + b_8(x)q = g_2(x) \end{cases} \quad (1)$$

where $a_j(x)$, $b_j(x)$, $g_1(x)$, $g_2(x)$ are given functions, and $a_j(x)$, $b_j(x)$ are continuous for $j = 1, 2, \ldots, 8$.

Let's have a look into system of third-order linear ordinary boundary value problems in one dimension for the pair of functions $p(x)$ and $q(x)$ in the following form [15]

$$\begin{cases} p''' + a_1(x)p'' + a_2(x)p' + a_3(x)p + a_4(x)q'' + a_5(x)q' + a_6(x)q = f(x), & a \leq x \leq b \\ q''' + b_1(x)q'' + b_2(x)q' + b_3(x)q + b_4(x)p'' + b_5(x)p' + b_6(x)p = g(x), & a \leq x \leq b. \end{cases} \quad (2)$$

Here, both of the equations are of third-order, therefore three boundary conditions are needed for both $p(x)$ and $q(x)$. Let us assume the boundary conditions at the ends, i.e.

$$p(a) = \alpha_1, \; p'(a) = \gamma_1, \; p'(b) = \beta_1, \; q(a) = \alpha_2, \; q'(a) = \gamma_2, \; q'(b) = \beta_2. \quad (3)$$

The terms $a_i(x)$ and $b_i(x)$ are continuous functions for $i = 1, 2, \ldots, n$ and, $f(x)$ and $g(x)$ are non-homogeneous.

The trial solutions for the two functions $p(x)$ and $q(x)$ of system (2) can be written as

$$\begin{cases} \tilde{p}(x) = \sum_{i=1}^{n-1} a_i \phi_i(x), & n \geq 1 \\ \tilde{q}(x) = \sum_{i=1}^{n-1} b_i \phi_i(x), & n \geq 1 \end{cases} \quad (4)$$

where $a_i$, $b_i$ are unknown parameters and $\phi_i(x)$ are basis functions which must satisfy the boundary conditions (3).

Using these approximations, we can define two residual functions $\varepsilon_p(x)$ and $\varepsilon_q(x)$:

$$\begin{cases} \varepsilon_p(x) = \tilde{p}''' + a_1(x)\tilde{p}'' + a_2(x)\tilde{p}' + a_3(x)\tilde{p} + a_4(x)\tilde{q}'' + a_5(x)\tilde{q}' + a_6(x)\tilde{q} - f(x) \\ \varepsilon_q(x) = \tilde{q}''' + b_1(x)\tilde{q}'' + b_2(x)\tilde{q}' + b_3(x)\tilde{q} + b_4(x)\tilde{p}'' + b_5(x)\tilde{p}' + b_6(x)\tilde{p} - g(x). \end{cases} \quad (5)$$

Now applying the Galerkin method, we get weighted residual system of equations

$$\begin{cases} \int_a^b \varepsilon_p(x)\phi_i(x)dx = 0 \\ \int_a^b \varepsilon_q(x)\phi_i(x)dx = 0. \end{cases} \quad (6)$$

Substituting $\varepsilon_p(x)$ and $\varepsilon_q(x)$, we get

$$\begin{cases} \int_a^b \left(\tilde{p}''' + a_1(x)\tilde{p}'' + a_2(x)\tilde{p}' + a_3(x)\tilde{p} + a_4(x)\tilde{q}'' + a_5(x)\tilde{q}' + a_6(x)\tilde{q}\right)\phi_i(x)dx = \int_a^b f(x)\phi_i(x)dx \\ \int_a^b (\tilde{q}''' + b_1(x)\tilde{q}'' + b_2(x)\tilde{q}' + b_3(x)\tilde{q} + b_4(x)\tilde{p}'' + b_5(x)\tilde{p}' + b_6(x)\tilde{p})\phi_i(x)dx = \int_a^b g(x)\phi_i(x)dx. \end{cases} \quad (7)$$

Applying integration by parts in the first term of (7) and setting $\phi_i(x) = 0$ at the boundary $x = a$ and $x = b$, we obtain

$$\int_a^b \frac{d^3 \tilde{p}}{dx^3} \phi_i(x)dx = \left[\frac{d\tilde{p}}{dx}\phi_i(x)\right]_a^b - \int_a^b \frac{d\tilde{p}}{dx}\frac{d\phi_i}{dx}dx$$

$$= \frac{d\tilde{p}}{dx}(x=b)\phi_i(b) - \frac{d\tilde{p}}{dx}(x=a)\phi_i(a) - \int_a^b \frac{d\tilde{p}}{dx}\frac{d\phi_i}{dx}dx$$

$$= -\left[\frac{d}{dx}\phi_i(x)\frac{d\tilde{p}}{dx}\right]_a^b + \int_a^b \frac{d^2}{dx^2}\phi_i(x)\frac{d\tilde{p}}{dx}dx$$

$$= -\left[\frac{d}{dx}\phi_i(x)\frac{d\tilde{p}}{dx}\right]_{x=b} + \left[\frac{d}{dx}\phi_i(x)\frac{d\tilde{p}}{dx}\right]_{x=a} + \int_a^b \frac{d^2}{dx^2}\phi_i(x)\frac{d\tilde{p}}{dx}dx$$

$$\int_a^b \frac{d^3 \tilde{q}}{dx^3}\phi_i(x)dx = -\left[\frac{d}{dx}\phi_i(x)\frac{d\tilde{q}}{dx}\right]_{x=b} + \left[\frac{d}{dx}\phi_i(x)\frac{d\tilde{q}}{dx}\right]_{x=a} + \int_a^b \frac{d^2}{dx^2}\phi_i(x)\frac{d\tilde{q}}{dx}dx.$$

Substituting these in (7), the system of weighted simultaneous residual equations is obtained

$$-\left[\frac{d\phi_i}{dx}\frac{d\tilde{p}}{dx}\right]_{x=b} + \left[\frac{d\phi_i}{dx}\frac{d\tilde{p}}{dx}\right]_{x=a} + \int_a^b \phi_i''(x)\tilde{p}'(x)dx$$

$$+ \int_a^b \left(a_1(x)\tilde{p}'' + a_2(x)\tilde{p}' + a_3(x)\tilde{p} + a_4(x)\tilde{q}'' + a_5(x)\tilde{q}' + a_6(x)\tilde{q}\right)\phi_i(x)dx = \int_a^b f(x)\phi_i(x)dx \quad (8a)$$

$$-\left[\frac{d\phi_i}{dx}\frac{d\tilde{q}}{dx}\right]_{x=b} + \left[\frac{d\phi_i}{dx}\frac{d\tilde{q}}{dx}\right]_{x=a} + \int_a^b \phi_i''(x)\tilde{q}'(x)dx$$

$$+ \int_a^b \left(b_1(x)\tilde{q}'' + b_2(x)\tilde{q}' + b_3(x)\tilde{q} + b_4(x)\tilde{p}'' + b_5(x)\tilde{p}' + b_6(x)\tilde{p}\right)\phi_i(x)dx = \int_a^b g(x)\phi_i(x)dx \quad (8b)$$

Applying the boundary conditions (3) in (8a) and (8b), we have respectively

$$-\left[\frac{d\phi_i}{dx}\beta_1\right]_{x=b} + \left[\frac{d\phi_i}{dx}\gamma_1\right]_{x=a} + \int_a^b \phi_i''(x)\tilde{p}'(x)dx$$
$$+ \int_a^b \left(a_1(x)\tilde{p}'' + a_2(x)\tilde{p}' + a_3(x)\tilde{p} + a_4(x)\tilde{q}'' + a_5(x)\tilde{q}' + a_6(x)\tilde{q}\right)\phi_i(x)dx = \int_a^b f(x)\phi_i(x)dx \quad (9a)$$

$$-\left[\frac{d\phi_i}{dx}\beta_2\right]_{x=b} + \left[\frac{d\phi_i}{dx}\gamma_2\right]_{x=a} + \int_a^b \phi_i''(x)\tilde{q}'(x)dx$$
$$+ \int_a^b \left(b_1(x)\tilde{q}'' + b_2(x)\tilde{q}' + b_3(x)\tilde{q} + b_4(x)\tilde{p}'' + b_5(x)\tilde{p}' + b_6(x)\tilde{p}\right)\phi_i(x)dx = \int_a^b g(x)\phi_i(x)dx \quad (9b)$$

Now putting trial solutions (4) into (9a) and (9b), and simplifying these for $i = 1, 2, \ldots, n$, we obtain

$$\sum_{j=1}^{n-1}\left(\int_a^b \left[\phi_i''(x)\phi_j'(x) + (a_1(x)\phi_j''(x) + a_2(x)\phi_j'(x) + a_3(x)\phi_j(x))\phi_i(x)\right]dx\right)a_j$$
$$+ \sum_{j=1}^{n-1}\left(\int_a^b (a_4(x)\phi_j''(x) + a_5(x)\phi_j'(x) + a_6(x)\phi_j(x))\phi_i(x)dx\right)b_j$$
$$= \int_a^b f(x)\phi_i(x)dx + \left[\frac{d\phi_i}{dx}\beta_1\right]_{x=b} - \left[\frac{d\phi_i}{dx}\gamma_1\right]_{x=a}, \quad (10a)$$

$$\sum_{j=1}^{n-1}\left(\int_a^b \left[\phi_i''(x)\phi_j'(x) + (b_1(x)\phi_j''(x) + b_2(x)\phi_j'(x) + b_3(x)\phi_j(x))\phi_i(x)\right]dx\right)b_j$$
$$+ \sum_{j=1}^{n-1}\left(\int_a^b (b_4(x)\phi_j''(x) + b_5(x)\phi_j'(x) + b_6(x)\phi_j(x))\phi_i(x)dx\right)a_j$$
$$= \int_a^b g(x)\phi_i(x)dx + \left[\frac{d\phi_i}{dx}\beta_2\right]_{x=b} - \left[\frac{d\phi_i}{dx}\gamma_2\right]_{x=a}. \quad (10b)$$

The matrix form is equivalent to the previously mentioned equations

$$\begin{cases} \sum_{j=1}^{n-1}(A_{j,i}a_j + H_{j,i}b_j) = F_i \\ \sum_{j=1}^{n-1}(C_{j,i}b_j + D_{j,i}a_j) = G_i \end{cases} \quad (11)$$

where,

$$A_{j,i} = \int_a^b \left[\phi_i''(x)\phi_j'(x) + (a_1(x)\phi_j''(x) + a_2(x)\phi_j'(x) + a_3(x)\phi_j(x))\phi_i(x)\right]dx$$

$$H_{j,i} = \int_a^b (a_4(x)\phi_j''(x) + a_5(x)\phi_j'(x) + a_6(x)\phi_j(x))\phi_i(x)dx$$

$$F_i = \int_a^b f(x)\phi_i(x)dx + \left[\frac{d\phi_i}{dx}\beta_1\right]_{x=b} - \left[\frac{d\phi_i}{dx}\gamma_1\right]_{x=a}$$

$$C_{j,i} = \int_a^b \left[\phi_i''(x)\phi_j'(x) + (b_1(x)\phi_j''(x) + b_2(x)\phi_j'(x) + b_3(x)\phi_j(x))\phi_i(x)\right] dx$$

$$D_{j,i} = \int_a^b (b_4(x)\phi_j''(x) + b_5(x)\phi_j'(x) + b_6(x)\phi_j(x))\phi_i(x) dx$$

$$G_i = g(x)\phi_i(x)dx + \left[\frac{d\phi_i}{dx}\beta_2\right]_{x=b} - \left[\frac{d\phi_i}{dx}\gamma_2\right]_{x=a}$$

$$i, j = 1, 2, \ldots, n$$

Now, we are going to solve nonlinear system of two third order differential equations in two unknown functions $p(x)$ and $q(x)$ of the form below is taking into consideration [8]

$$\begin{cases} a_1(x)p''' + a_2(x)p'' + a_3(x)p' + a_4(x)p + a_5(x)q''' + a_6(x)q'' + a_7(x)q' + a_8(x)q + M_1(p,q) = g_1(x) \\ b_1(x)p''' + b_2(x)p'' + b_3(x)p' + b_4(x)p + b_5(x)q''' + b_6(x)q'' + b_7(x)q' + b_8(x)q + M_2(p,q) = g_2(x) \end{cases} \quad (12)$$

where $a_j(x)$, $b_j(x)$, $g_1(x)$, $g_2(x)$ are given functions, $M_1$, $M_2$ represent nonlinear functions and $a_j(x)$, $b_j(x)$ are continuous for $j = 1, 2, \ldots, 8$. Since the system of third-order boundary value problems consist of nonlinear terms, then we can't solve the system directly. In this case, mathematical formulation depends on the problem. In order to get the initial values of the parameters, we must neglect the nonlinear terms and solve the linear system. After getting the initial values of the parameters we will be able to solve the system. Then putting the parameters into the trial solutions, we will obtain the approximate solutions for the functions $p(x)$ and $q(x)$. The details are described in the following section.

## 4 Numerical Results & Discussions

In this study, we use four systems; one linear and three nonlinear, which are available in the existing literature. To verify the effectiveness of the derived formulations, Dirichlet and Neumann boundary conditions are considered. For each case we find the approximate solutions using different number of parameters with Bernstein polynomials, and we compare these solutions with the exact solutions, and graphically which are shown in the same figures. Since the convergence of linear BVP is calculated by

$$E = |\tilde{p}_{n+1}(x) - \tilde{p}_n(x)| < \delta_1$$

where $\tilde{p}_n(x)$ denotes the approximate solution using $n$ polynomials and $\delta_1$ depends on the problem. In this case, $\delta_1 < 10^{-8}$. In addition, the convergence of nonlinear BVP is assumed when the absolute error of two consecutive iterations, $\delta_2$ satisfies

$$E = |\tilde{p}_n^{N+1}(x) - \tilde{p}_n^N(x)| < \delta_2$$

and, in this case, $\delta_2 < 0^{-10}$.

**Example 1.** Consider the following system of third-order nonlinear boundary value problem [21, 22]

$$\begin{cases} p'''(x) + 2p'(x) + xq(x) = f(x), & 0 < x < 1 \\ q'''(x) + \dfrac{p''(x)q''(x)}{6} = g(x), & 0 < x < 1 \\ p(0) = q(0) = 0, \ p(1) = q(1) = 0, \ p'(0) = q'(0) = 0 \end{cases} \quad (13)$$

where

$$f(x) = x^5 - x^3 - 18x^2 + 12x - 18$$

$$g(x) = -36x^3 + 12x^2 + 30x - 2$$

The exact solutions are $p(x) = 3x^2 - 3x^3$ and $q(x) = x^4 - x^2$.

Here, we use Bernstein polynomials as trial approximate solution to solve the system (13). Let us consider the trial solution of the form

$$\begin{cases} \tilde{p}(x) = \theta_0(x) + \sum_{i=1}^{n-1} a_i \phi_i(x) \\ \tilde{q}(x) = \phi_0(x) + \sum_{i=1}^{n-1} b_i \phi_i(x) \end{cases} \quad (14)$$

where $a_i$ and $b_i$ are parameters and $\phi_i(x)$ are trial functions (Bernstein polynomials) which satisfy the boundary conditions. Here, we can choose $\theta_0(x) = 0$ and $\phi_0(x) = 0$ because the boundary conditions are homogeneous.

Now applying the Galerkin method, we get weighted residual system of equations

$$\begin{cases} \int_0^1 \left( \tilde{p}'''(x) + 2\tilde{p}'(x) + x\tilde{q}(x) \right) \phi_i(x) dx = \int_0^1 f(x) \phi_i(x) dx \\ \int_0^1 \left( \tilde{q}'''(x) + \dfrac{\tilde{p}''(x)\tilde{q}''(x)}{6} \right) \phi_i(x) dx = \int_0^1 g(x) \phi_i(x) dx \end{cases} \quad (15)$$

Applying integration by parts in the first term of (22), we obtain

$$\int_0^1 \frac{d^3 \tilde{p}}{dx^3} \phi_i(x) dx = \left[ \frac{d\tilde{p}}{dx} \phi_i(x) \right]_0^1 - \int_0^1 \frac{d\tilde{p}}{dx} \frac{d\phi_i}{dx} dx$$

$$= \frac{d\tilde{p}}{dx}(x=1)\phi_i(1) - \frac{d\tilde{p}}{dx}(x=0)\phi_i(0) - \int_0^1 \frac{d\tilde{p}}{dx} \frac{d\phi_i}{dx} dx$$

$$= -\left[ \frac{d}{dx}\phi_i(x) \frac{d\tilde{p}}{dx} \right]_0^1 + \int_0^1 \frac{d^2}{dx^2} \phi_i(x) \frac{d\tilde{p}}{dx} dx$$

$$\int_0^1 \frac{d^3 \tilde{q}}{dx^3} \phi_i(x) dx = -\left[ \frac{d}{dx}\phi_i(x) \frac{d\tilde{q}}{dx} \right]_0^1 + \int_0^1 \frac{d^2}{dx^2} \phi_i(x) \frac{d\tilde{q}}{dx} dx$$

This system can be converted into modified Galerkin form in the usual sense, and using (14) to obtain

$$\sum_{j=1}^{n-1} \left( \int_0^1 \left[ \phi_j'(x)\phi_i''(x) + 2\phi_j'(x)\phi_i(x) \right] dx - [\phi_j'(x)\phi_i'(x)]_{x=1} \right) a_j$$

$$+ \sum_{j=1}^{n-1} \left( \int_0^1 x\phi_j(x)\phi_i(x) dx \right) b_j = \int_0^1 f(x)\phi_i(x) dx, \quad i = 1, 2, \ldots, n \quad (16a)$$

$$\sum_{j=1}^{n-1} \left( \int_0^1 \phi_j'(x)\phi_i''(x)dx - [\phi_j'(x)\phi_i'(x)]_{x=1} \right) b_j$$

$$+ \frac{1}{6}\int_0^1 \left(\sum_{k=1}^{n-1} a_k \phi_k''(x)\right) \left(\sum_{k=1}^{n-1} b_k \phi_k''(x)\right) \phi_i(x)dx = \int_0^1 g(x)\phi_i(x)dx, \quad i = 1,2,\ldots,n \quad (16b)$$

The above equations are equivalent to the matrix form

$$\sum_{j=1}^{n-1} (A_{j,i}a_j + H_{j,i}b_j) = F_i, \quad i = 1,2,\ldots,n \quad (17)$$

$$\sum_{j=1}^{n-1} C_{j,i}b_j = G_i + N_i, \quad i = 1,2,\ldots,n \quad (18)$$

in which

$$A_{j,i} = \int_0^1 \left( \phi_j'(x)\phi_i''(x) + 2\phi_j'(x)\phi_i(x) \right) dx - [\phi_j'(x)\phi_i'(x)]_{x=1}$$

$$H_{j,i} = \int_0^1 x\phi_j(x)\phi_i(x)dx$$

$$F_i = \int_0^1 (x^5 - x^3 - 18x^2 + 12x - 18)\phi_i(x)dx$$

$$C_{j,i} = \int_0^1 \phi_j'(x)\phi_i''(x)dx - [\phi_j'(x)\phi_i'(x)]_{x=1}$$

$$G_i = \int_0^1 (-36x^3 + 12x^2 + 30x - 2)\phi_i(x)dx$$

$$N_i = -\frac{1}{6}\int_0^1 \left(\sum_{k=1}^n a_k \phi_k''(x)\right) \left(\sum_{k=1}^n b_k \phi_k''(x)\right) \phi_i(x)dx$$

$$i, j = 1,2,\ldots,n$$

Neglecting the nonlinear term $N_i$ in (18), the initial values of $a_j$ and $b_j$ are obtained. Therefore, we obtain initial values solving the system

$$\begin{cases} \sum_{j=1}^{n-1} (A_{j,i}a_j + H_{j,i}b_j) = F_i, & i = 1,2,\ldots,n \\ \sum_{j=1}^{n-1} C_{j,i}b_j = G_i & i = 1,2,\ldots,n \end{cases} \quad (19)$$

Then we substitute these values into (17) and (18) and obtain the new values of $a_j$ and $b_j$. The iterative process keeps going till the converged estimates of the unknown coefficients are achieved. We have an approximation to the BVP (13) by replacing the final quantities in (14).

Satisfying the homogeneous boundary conditions and using degree of polynomial 3, we may obtain the approximate solutions of $p(x)$ and $q(x)$ are

$$\tilde{p}(x) = 0.00054548x(1-x)^2 + 2.99843577x^2(1-x)$$

and

$$\tilde{q}(x) = 0.39311569x(1-x)^2 - 2.07669616x^2(1-x),$$

respectively.

The following Table 1 shows the numerical outcomes for the given problem. Here we have used 2, 3 and 4 Bernstein polynomials in column 2,3 and 4 respectively. Consider the fact that we iterate four times in order to achieve the approximations.

Table 1: Absolute errors $|p(x) - \tilde{p}(x)|$ for **Example 1**.

| $x$ | Present Method (GWRM) | | | HPM [21] |
|---|---|---|---|---|
| | **Degree 3** | **Degree 4** | **Degree 5** | |
| 0.1 | $3.010565 \times 10^{-5}$ | $8.406867 \times 10^{-10}$ | $1.680418 \times 10^{-9}$ | $8.74250 \times 10^{-8}$ |
| 0.2 | $1.976596 \times 10^{-5}$ | $3.596601 \times 10^{-8}$ | $2.284382 \times 10^{-9}$ | $3.41082 \times 10^{-5}$ |
| 0.3 | $1.836086 \times 10^{-5}$ | $9.477548 \times 10^{-8}$ | $5.387116 \times 10^{-9}$ | $7.35812 \times 10^{-5}$ |
| 0.4 | $7.161663 \times 10^{-5}$ | $1.653393 \times 10^{-7}$ | $1.220891 \times 10^{-8}$ | $1.23210 \times 10^{-4}$ |
| 0.5 | $1.273431 \times 10^{-4}$ | $2.343983 \times 10^{-7}$ | $2.200062 \times 10^{-8}$ | $6.24822 \times 10^{-1}$ |
| 0.6 | $1.728822 \times 10^{-4}$ | $2.873641 \times 10^{-7}$ | $3.242916 \times 10^{-8}$ | $2.31537 \times 10^{-4}$ |
| 0.7 | $1.955756 \times 10^{-4}$ | $3.083190 \times 10^{-7}$ | $3.996300 \times 10^{-8}$ | $2.75507 \times 10^{-5}$ |
| 0.8 | $1.827652 \times 10^{-4}$ | $2.800157 \times 10^{-7}$ | $4.025768 \times 10^{-8}$ | $2.92801 \times 10^{-4}$ |
| 0.9 | $1.217927 \times 10^{-4}$ | $1.838780 \times 10^{-7}$ | $2.854130 \times 10^{-8}$ | $2.39684 \times 10^{-4}$ |

Here, we see that the maximum absolute error of $p(x)$ using 3, 4 and 5 degree polynomials are $1.955756 \times 10^{-4}$, $3.083190 \times 10^{-7}$ and $4.025768 \times 10^{-8}$, respectively whereas the maximum absolute errors using Homotopy Perturbation Method (HPM) in [21] and Block method in [22] are $6.24822 \times 10^{-1}$ and $6.25 \times 10^{-4}$, respectively.

Table 2: Absolute errors $|q(x) - \tilde{q}(x)|$ for **Example 1**.

| $x$ | Present Method (GWRM) | | | HPM [21] |
|---|---|---|---|---|
| | **Degree 3** | **Degree 4** | **Degree 5** | |
| 0.1 | $2.305211 \times 10^{-2}$ | $1.390676 \times 10^{-8}$ | $3.551924 \times 10^{-7}$ | $1.34169 \times 10^{-5}$ |
| 0.2 | $2.226453 \times 10^{-2}$ | $7.162991 \times 10^{-6}$ | $1.658188 \times 10^{-6}$ | $5.22346 \times 10^{-5}$ |
| 0.3 | $8.856149 \times 10^{-3}$ | $1.724459 \times 10^{-5}$ | $2.285279 \times 10^{-6}$ | $1.12022 \times 10^{-4}$ |
| 0.4 | $8.354171 \times 10^{-3}$ | $2.692795 \times 10^{-5}$ | $4.263502 \times 10^{-7}$ | $9.61854 \times 10^{-2}$ |
| 0.5 | $2.294756 \times 10^{-2}$ | $3.375417 \times 10^{-5}$ | $3.196644 \times 10^{-6}$ | $2.00026 \times 10^{-1}$ |
| 0.6 | $3.090514 \times 10^{-2}$ | $3.613628 \times 10^{-5}$ | $6.873504 \times 10^{-6}$ | $3.30775 \times 10^{-4}$ |
| 0.7 | $3.060805 \times 10^{-2}$ | $3.335918 \times 10^{-5}$ | $8.739233 \times 10^{-6}$ | $3.78034 \times 10^{-4}$ |
| 0.8 | $2.283741 \times 10^{-2}$ | $2.557966 \times 10^{-5}$ | $7.607479 \times 10^{-6}$ | $3.86227 \times 10^{-4}$ |
| 0.9 | $1.077435 \times 10^{-2}$ | $1.382641 \times 10^{-5}$ | $3.803978 \times 10^{-6}$ | $3.11011 \times 10^{-4}$ |

On the other hand, Table 2 shows that the maximum absolute error of $q(x)$ using 3, 4 and 5 degree polynomials are $3.090514 \times 10^{-2}$, $3.613628 \times 10^{-5}$ and $8.739233 \times 10^{-6}$, respectively. The maximum absolute errors using HPM in [21] and Block method in [22] are $2.00026 \times 10^{-1}$ and $1.48 \times 10^{-2}$, respectively. Thus, our proposed method reflects the better results than the previous results.

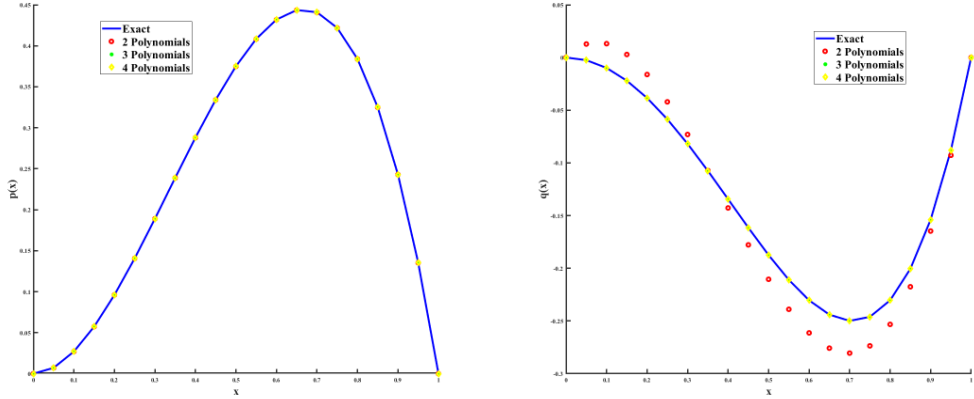

Figure 1: Exact and approximate solutions of $p(x)$ and $q(x)$ for Example 1.

**Example 2.** Consider the following nonlinear system of third-order boundary value problem [21, 22]

$$\begin{cases} p'''(x) - 4q''(x) + p''(x)q'(x) = f(x), & 0 < x < 1 \\ q'''(x) + 4q'(x) - p''(x) + p'(x)q''(x) = g(x), & 0 < x < 1 \\ p(0) = q(0) = 0, \ p(1) = q(1) = 1, \ p'(0) = q'(0) = 0 \end{cases} \quad (20)$$

where

$$f(x) = 36x^4 \quad \text{and} \quad g(x) = 24x^4 + 6$$

The exact solutions are $p(x) = x^4$ and $q(x) = x^3$.

Let us consider the trial solution of the form

$$\begin{cases} \tilde{p}(x) = \theta_0(x) + \sum_{i=1}^{n} a_i \phi_i(x) \\ \tilde{q}(x) = \phi_0(x) + \sum_{i=1}^{n} b_i \phi_i(x) \end{cases} \quad (21)$$

where $a_i$ and $b_i$ are parameters and $\phi_i(x)$ are trial functions (Bernstein polynomials) which satisfy the boundary conditions. Here, we can choose $\theta_0(x) = x^2$ and $\phi_0(x) = x^2$ in order to satisfy the boundary conditions.

Now applying the Galerkin method, we get weighted residual system of equations

$$\begin{cases} \int_0^1 \left( \tilde{p}'''(x) - 4\tilde{q}''(x) + \tilde{p}''(x)\tilde{q}'(x) \right) \phi_i(x) dx = \int_0^1 f(x)\phi_i(x) dx \\ \int_0^1 \left( \tilde{q}'''(x) + 4\tilde{q}'(x) - \tilde{p}''(x) + \tilde{p}'(x)\tilde{q}''(x) \right) \phi_i(x) dx = \int_0^1 g(x)\phi_i(x) dx \end{cases} \quad (22)$$

This system can be converted to modified Galerkin form in the usual way and using (21) we finally obtain the matrix form

$$\sum_{j=1}^{n} (A_{j,i} a_j + H_{j,i} b_j) = F_i + N_i, \quad i = 1, 2, \ldots, n \quad (23)$$

$$\sum_{j=1}^{n} (C_{j,i} b_j + D_{j,i} a_j) = G_i + M_i, \quad i = 1, 2, \ldots, n \quad (24)$$

where

$$A_{j,i} = \int_0^1 \left[\phi_j'(x)\phi_i''(x) + 2x\phi_j''(x)\phi_i(x)\right] dx - [\phi_j'(x)\phi_i'(x)]_{x=1}$$

$$H_{j,i} = \int_0^1 (4\phi_j'(x)\phi_i'(x) + 2\phi_j'(x)\phi_i(x))dx$$

$$F_i = \int_0^1 \left(f(x)\phi_i(x) - 4x\phi_i(x) - 8x\phi_i'(x) - 2x\phi_i''(x)\right) dx + [2x\phi_i'(x)]_{x=1}$$

$$N_i = -\int_0^1 \left(\sum_{k=1}^n a_k\phi_k''(x)\right)\left(\sum_{k=1}^n b_k\phi_k'(x)\right)\phi_i(x)dx$$

$$C_{j,i} = \int_0^1 (\phi_j'(x)\phi_i''(x) + 4\phi_j'(x)\phi_i(x) + 2x\phi_j''(x)\phi_i(x))dx - [\phi_j'(x)\phi_i'(x)]_{x=1}$$

$$D_{j,i} = \int_0^1 (4\phi_j'(x)\phi_i'(x) + 2\phi_j'(x)\phi_i(x))dx$$

$$G_i = \int_0^1 \left(g(x)\phi_i(x) - 12x\phi_i(x) - 2x\phi_i'(x) - 2x\phi_i''(x)\right) dx + [2x\phi_i'(x)]_{x=1}$$

$$M_i = -\int_0^1 \left(\sum_{k=1}^n a_k\phi_k'(x)\right)\left(\sum_{k=1}^n b_k\phi_k''(x)\right)\phi_i(x)dx$$

$$i,j = 1,2,\ldots,n$$

Neglecting the nonlinear terms $N_i$ and $M_i$ in (23) and (24) respectively, the initial values of $a_j$ and $b_j$ are obtained. Therefore, we obtain initial values solving the system

$$\begin{cases} \sum_{j=1}^n (A_{j,i}a_j + H_{j,i}b_j) = F_i, & i = 1,2,\ldots,n \\ \sum_{j=1}^n (C_{j,i}b_j + D_{j,i}a_j) = G_i & i = 1,2,\ldots,n \end{cases} \quad (25)$$

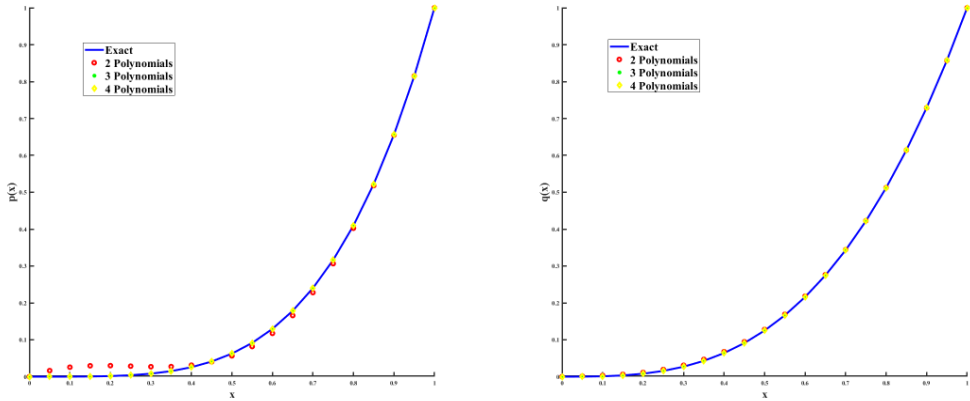

Figure 2: Exact and approximate solutions of $p(x)$ and $q(x)$ for Example 2.

Then we substitute these values in (23) and (24) and obtain the new values of $a_j$ and $b_j$. The iterative process keeps going till the converged estimates of the unknown coefficients are achieved. We have an approximation to the BVP (20) by replacing the final quantities in (21). The approximate solutions of $p(x)$ and $q(x)$ using two parameters with polynomial degree 3 are:

$$\tilde{p}(x) = 2.36893939x^3 - 1.77969288x^2 + 0.41075348x$$

and
$$\tilde{q}(x) = 1.03113495x^3 - 0.05887466x^2 + 0.02773971x,$$
respectively.

Table 3 and Table 4 show the numerical outcomes for the given problem. Here we have used 2, 3 and 4 Bernstein polynomials in column 2, 3 and 4 respectively. Consider the fact that we iterate four times in order to achieve the approximations. Here, we see that the maximum absolute error of $p(x)$ in Table 3 using 2, 3 and 4 polynomials are $2.831450 \times 10^{-2}$, $1.616954 \times 10^{-7}$ and $8.469325 \times 10^{-8}$, respectively whereas the maximum absolute errors using Homotopy Perturbation Method (HPM) in [21] and Block method in [22] are $7.66527 \times 10^{-2}$ and $5.50 \times 10^{-3}$, respectively.

Table 3: Absolute errors for $|p(x) - \tilde{p}(x)|$ for **Example 2**.

| $x$ | Present Method (GWRM) | | | HPM [21] |
|---|---|---|---|---|
| | **Degree 3** | **Degree 4** | **Degree 5** | |
| 0.1 | $2.554736 \times 10^{-2}$ | $1.060524 \times 10^{-7}$ | $7.896618 \times 10^{-9}$ | $8.74250 \times 10^{-8}$ |
| 0.2 | $2.831450 \times 10^{-2}$ | $6.744158 \times 10^{-8}$ | $1.381743 \times 10^{-8}$ | $1.08399 \times 10^{-3}$ |
| 0.3 | $1.891505 \times 10^{-2}$ | $2.037354 \times 10^{-8}$ | $2.377000 \times 10^{-8}$ | $4.71442 \times 10^{-3}$ |
| 0.4 | $5.562655 \times 10^{-3}$ | $9.053284 \times 10^{-8}$ | $3.918424 \times 10^{-8}$ | $1.14828 \times 10^{-2}$ |
| 0.5 | $5.929053 \times 10^{-3}$ | $1.047750 \times 10^{-7}$ | $5.790583 \times 10^{-8}$ | $2.18865 \times 10^{-2}$ |
| 0.6 | $1.214644 \times 10^{-2}$ | $5.343750 \times 10^{-8}$ | $7.518965 \times 10^{-8}$ | $3.5986 \times 10^{-2}$ |
| 0.7 | $1.207586 \times 10^{-2}$ | $4.454330 \times 10^{-8}$ | $8.469325 \times 10^{-8}$ | $6.87692 \times 10^{-2}$ |
| 0.8 | $7.103685 \times 10^{-3}$ | $1.416323 \times 10^{-7}$ | $7.947025 \times 10^{-8}$ | $7.66527 \times 10^{-2}$ |
| 0.9 | $1.016277 \times 10^{-3}$ | $1.616954 \times 10^{-7}$ | $5.296380 \times 10^{-8}$ | $6.20685 \times 10^{-2}$ |

On the other hand, Table 4 shows that the maximum absolute error of $q(x)$ using 2, 3 and 4 polynomials are $3.863836 \times 10^{-3}$, $2.553118 \times 10^{-7}$ and $3.153957 \times 10^{-8}$, respectively. The maximum absolute errors using HPM in [21] and Block method in [22] are $6.57775 \times 10^{-2}$ and $2.39 \times 10^{-2}$, respectively. In this example, we show that the proposed method is far better than the existing methods in [21] and in [22].

Table 4: Absolute errors for $|q(x) - \tilde{q}(x)|$ for **Example 2**.

| $x$ | Present Method (GWRM) | | | HPM [21] |
|---|---|---|---|---|
| | **Degree 3** | **Degree 4** | **Degree 5** | |
| 0.1 | $2.216359 \times 10^{-3}$ | $8.201831 \times 10^{-8}$ | $5.335278 \times 10^{-9}$ | $6.73806 \times 10^{-4}$ |
| 0.2 | $3.442034 \times 10^{-3}$ | $7.169045 \times 10^{-8}$ | $1.268373 \times 10^{-9}$ | $2.82326 \times 10^{-3}$ |
| 0.3 | $3.863836 \times 10^{-3}$ | $3.608265 \times 10^{-8}$ | $7.607422 \times 10^{-9}$ | $6.61192 \times 10^{-3}$ |
| 0.4 | $3.668574 \times 10^{-3}$ | $1.982045 \times 10^{-8}$ | $1.026833 \times 10^{-8}$ | $1.21045 \times 10^{-3}$ |
| 0.5 | $3.043057 \times 10^{-3}$ | $4.508859 \times 10^{-8}$ | $1.119354 \times 10^{-8}$ | $1.91209 \times 10^{-2}$ |
| 0.6 | $2.174096 \times 10^{-3}$ | $1.116311 \times 10^{-7}$ | $1.424742 \times 10^{-8}$ | $2.69597 \times 10^{-2}$ |
| 0.7 | $1.248500 \times 10^{-3}$ | $1.967514 \times 10^{-7}$ | $2.178228 \times 10^{-8}$ | $3.39306 \times 10^{-2}$ |
| 0.8 | $4.530788 \times 10^{-4}$ | $2.553118 \times 10^{-7}$ | $3.120430 \times 10^{-8}$ | $6.57775 \times 10^{-2}$ |
| 0.9 | $2.535774 \times 10^{-5}$ | $2.197343 \times 10^{-7}$ | $3.153957 \times 10^{-8}$ | $2.88236 \times 10^{-2}$ |

# 5 Application

Although in [24], Agarwal has extensively covered the theorems of the criteria for the existence and uniqueness of solutions of the sixth-order BVPs in a book, it does not include any numerical techniques. Islam and Hossain [23] solved the sixth-order BVPs using Galerkin method. Modified decomposition method was used in [25] to find the solution of the sixth-order BVPs by Wazwaz. Aasma Khalid et al. used Cubic B-spline in [26] and [27] in order to solve the linear and nonlinear sixth-order BVPs, respectively. Cubic-nonpolynomial spline (CNPS) and cubic-polynomial spline (CPS) were used to obtain the solutions of BVPs in [28] arising in hydrodynamic and magnetohydro-dynamic stability theory. Noor and Mohyud-Din [29] solved the sixth-order BVP using homotopy perturbation method. However, in this section, we show that higher even order BVP may be solved, in the alternative way, by the method of reduction order into system of lower order BVPs. For this, we experiment the proposed method to sixth order BVP.

**Example 3.** Consider the linear sixth-order boundary value problem [23, 26, 29]

$$\frac{d^6 p}{dx^6} - p = -6e^x, \quad 0 \leq x \leq 1 \tag{26}$$

Subject to the boundary conditions

$$p(0) = 1, \ p(1) = 0, \ p'(0) = 0, \ p'(1) = -e, \ p''(0) = -1, \ p''(1) = -2e$$

The analytic solution of the above problem is, $p(x) = (1-x)e^x$.

If we introduce a new function, say, $q(x)$ such that

$$q = \frac{d^3 p}{dx^3} \tag{27}$$

then (30) is clearly equivalent to two third-order differential equations

$$\begin{cases} \frac{d^3 q}{dx^3} - p = -6e^x \\ \frac{d^3 p}{dx^3} - q = 0 \end{cases} \tag{28}$$

Table 5: Exact, approximate and absolute errors of **Example 3**.

| $x$ | Present Method (GWRM) | | | Cubic B-Spline [26] | HPM [29] |
|---|---|---|---|---|---|
| | **Exact** | **Approximate** | **Abs Error** | | |
| 0.1 | 0.99465383 | 0.99464299 | $1.08 \times 10^{-5}$ | $1.18 \times 10^{-5}$ | $4.09 \times 10^{-4}$ |
| 0.2 | 0.97712221 | 0.97713362 | $1.14 \times 10^{-5}$ | $4.29 \times 10^{-5}$ | $7.78 \times 10^{-4}$ |
| 0.3 | 0.94490117 | 0.94491960 | $1.84 \times 10^{-5}$ | $8.53 \times 10^{-5}$ | $1.07 \times 10^{-3}$ |
| 0.4 | 0.89509482 | 0.89510135 | $6.53 \times 10^{-6}$ | $1.28 \times 10^{-4}$ | $1.26 \times 10^{-3}$ |
| 0.5 | 0.82436064 | 0.82435116 | $9.48 \times 10^{-6}$ | $1.59 \times 10^{-4}$ | $1.32 \times 10^{-3}$ |
| 0.6 | 0.72884752 | 0.72883236 | $1.52 \times 10^{-5}$ | $1.67 \times 10^{-4}$ | $1.26 \times 10^{-3}$ |
| 0.7 | 0.60412581 | 0.60411856 | $7.26 \times 10^{-6}$ | $1.45 \times 10^{-4}$ | $1.07 \times 10^{-3}$ |
| 0.8 | 0.44510819 | 0.44511280 | $4.61 \times 10^{-6}$ | $9.47 \times 10^{-5}$ | $7.78 \times 10^{-4}$ |
| 0.9 | 0.24596031 | 0.24596678 | $6.46 \times 10^{-6}$ | $3.33 \times 10^{-5}$ | $4.09 \times 10^{-4}$ |

with boundary conditions

$$p(0) = 1, \ p(1) = 0, \ p'(0) = 0, \ q(0) = -2, \ q(1) = -3e, \ q'(0) = -3 \quad (29)$$

The following Table 5 shows the numerical outcomes for the given problem. We have used 4 Bernstein polynomials of degree 5. Here, we see that the maximum absolute error of $p(x)$ in Table 5 is $1.84 \times 10^{-5}$ whereas the maximum absolute errors using Cubic B-Spline method in [26] and Homotopy Perturbation Method (HPM) in [29] are $1.67 \times 10^{-4}$ and $1.32 \times 10^{-3}$, respectively.

Using the method illustrated in the previous section and for different number of polynomials, the maximum absolute errors and the previous results available in the literature are summarized in Table 6. The accuracy of the present method is remarkable.

Table 6: Absolute errors for $|p(x) - \tilde{p}(x)|$ for **Example 3**.

| Number of Polynomial used | Present Method (GWRM) | In [23] |
|---|---|---|
| 10 | $1.165 \times 10^{-14}$ | $1.126 \times 10^{-13}$ |
| 11 | $1.986 \times 10^{-16}$ | $2.311 \times 10^{-15}$ |
| 12 | $4.771 \times 10^{-18}$ | $2.220 \times 10^{-16}$ |
| 13 | $7.228 \times 10^{-20}$ | $2.220 \times 10^{-16}$ |

**Example 4.** Consider the nonlinear sixth-order boundary value problem [25, 28, 27, 29]

$$\frac{d^6 p}{dx^6} = e^{-x} p^2(x), \ 0 \leq x \leq 1 \quad (30)$$

Subject to the boundary conditions

$$p(0) = p''(0) = p^{iv}(0) = 1, \ p(1) = p''(1) = p^{iv}(1) = e$$

The analytic solution of the above problem is, $p(x) = e^x$.

Introducing a new function, $q(x)$, equation (30) is equivalent to two third-order differential equations

$$\begin{cases} \dfrac{d^3 q}{dx^3} = e^{-x} p^2 \\ \dfrac{d^3 p}{dx^3} - q = 0 \end{cases} \quad (31)$$

Table 7: Exact, approximate and absolute errors of **Example 4**.

| $x$ | Exact | Approximate | Abs Error | MDM [25] | Cubic B-Spline [27] | CPS [28] | HPM [29] |
|---|---|---|---|---|---|---|---|
| 0.1 | 1.10517092 | 1.10516092 | $1.00 \times 10^{-5}$ | $1.23 \times 10^{-4}$ | $3.95 \times 10^{-6}$ | $5.05 \times 10^{-5}$ | $1.23 \times 10^{-4}$ |
| 0.2 | 1.22140276 | 1.22141684 | $1.41 \times 10^{-5}$ | $2.35 \times 10^{-4}$ | $1.43 \times 10^{-5}$ | $9.23 \times 10^{-5}$ | $2.35 \times 10^{-4}$ |
| 0.3 | 1.34985881 | 1.34989851 | $3.97 \times 10^{-5}$ | $3.25 \times 10^{-4}$ | $2.79 \times 10^{-5}$ | $1.25 \times 10^{-4}$ | $3.25 \times 10^{-4}$ |
| 0.4 | 1.49182470 | 1.49187513 | $5.04 \times 10^{-5}$ | $3.85 \times 10^{-4}$ | $4.07 \times 10^{-5}$ | $1.47 \times 10^{-4}$ | $3.85 \times 10^{-4}$ |
| 0.5 | 1.64872127 | 1.64876483 | $4.35 \times 10^{-5}$ | $4.08 \times 10^{-4}$ | $4.88 \times 10^{-5}$ | $1.59 \times 10^{-4}$ | $4.08 \times 10^{-4}$ |
| 0.6 | 1.82211880 | 1.82214507 | $2.63 \times 10^{-5}$ | $3.91 \times 10^{-4}$ | $4.92 \times 10^{-5}$ | $1.58 \times 10^{-4}$ | $3.92 \times 10^{-4}$ |
| 0.7 | 2.01375271 | 2.01376310 | $1.04 \times 10^{-5}$ | $3.36 \times 10^{-4}$ | $4.09 \times 10^{-5}$ | $1.44 \times 10^{-4}$ | $3.36 \times 10^{-4}$ |
| 0.8 | 2.22554093 | 2.22554635 | $5.42 \times 10^{-6}$ | $2.45 \times 10^{-4}$ | $2.56 \times 10^{-5}$ | $1.14 \times 10^{-4}$ | $2.46 \times 10^{-4}$ |
| 0.9 | 2.45960311 | 2.45961289 | $9.77 \times 10^{-6}$ | $1.29 \times 10^{-4}$ | $8.63 \times 10^{-6}$ | $6.70 \times 10^{-5}$ | $1.29 \times 10^{-4}$ |

Table 7 shows the numerical outcomes for the given problem. We have used 4 Bernstein polynomials of degree 5. Here, we see that the maximum absolute error of $p(x)$ in Table 7 is $5.04 \times 10^{-5}$ whereas the maximum absolute errors using Modified Decomposition Method (MPM) in [25], Cubic B-Spline method in [27], Cubic Polynomial Spline (CPS) in [28] and Homotopy Perturbation Method in [29] are $4.08 \times 10^{-4}$, $4.92 \times 10^{-5}$, $1.59 \times 10^{-4}$ and $4.08 \times 10^{-4}$, respectively. This concludes that the present method may be exploited with considerable accuracy.

# 6 Conclusion

We have deduced the formulation of the Galerkin weighted residual method for system of third-order boundary value problems in details. We can determine the solutions at each point within the problem's domain by using this method. Bernstein polynomials have been used in this method as the trial functions in the approximation. The focus has been on the formulations as well as on the performance of the accuracy. Some examples are tested to verify the effectiveness of the desired formulations. A good agreement has been established when comparing the approximate solutions with the exact solutions. We have shown that this method may be applied to higher-order systems and/or single BVPs to get the desired accuracy.